# Homogenization of a singular random one-dimensional PDE

Bogdan Iftimie[a], Étienne Pardoux[b] and Andrey Piatnitski[c]

[a]*Department of Mathematics, Academy of Economical Studies, 70167 Bucharest, Romania. E-mail: iftimieb@csie.ase.ro*
[b]*LATP, CMI, Universite de Provence, 39 rue Joliot-Curie, 13453 Marseille, Cedex 13013, France.
E-mail: Etienne.Pardoux@cmi.univ-mrs.fr*
[c]*Narvik University College, Postbox 385, 8505 Narvik, Norway and P.N. Lebedev Physical Institute RAS, Leninski prospect 53, 119991 Moscow, Russia. E-mail: andrey@sci.lebedev.ru*

**Abstract.** This paper deals with the homogenization problem for a one-dimensional parabolic PDE with random stationary mixing coefficients in the presence of a large zero order term. We show that under a proper choice of the scaling factor for the said zero order terms, the family of solutions of the studied problem converges in law, and describe the limit process. It should be noted that the limit dynamics remain random.

**Résumé.** Cet article traite de l'homogénéisation d'une équation aux dérivées partielles en dimension un d'espace, avec des coefficients aléatoires stationnaires et mélangeants, en présence d'u terme d'ordre zéro fortement oscillant. Nous montrons qu'avec un choix convenable du facteur d'échelle de ce terme d'ordre zéro, les solutions du problème étudié convergent en loi, et nous décrivons le processus limite. On peut noter que la dynamique limite est elle aussi aléatoire.

## 1. Introduction

Our goal is to study the limit, as $\varepsilon \to 0$, of a linear parabolic PDE of the form

$$\frac{\partial u^\varepsilon}{\partial t}(t,x) = \frac{1}{2}\frac{\partial}{\partial x}\left(a\left(\frac{\cdot}{\varepsilon}\right)\frac{\partial u^\varepsilon}{\partial x}\right)(t,x) + \frac{1}{\sqrt{\varepsilon}}c\left(\frac{x}{\varepsilon}\right)u^\varepsilon(t,x), \quad t \geq 0, x \in \mathbb{R};$$
$$u^\varepsilon(0,x) = g(x), \quad x \in \mathbb{R},$$
(1.1)

where $a$ and $c$ are stationary random fields, and $c$ is centered.

Let us recall (see [1]) that in the periodic case the equation

$$\frac{\partial}{\partial t}u = \frac{1}{2}\frac{\partial}{\partial x}\left(a\left(\frac{x}{\varepsilon}\right)\frac{\partial}{\partial x}u\right) + \frac{1}{\varepsilon}c\left(\frac{x}{\varepsilon}\right)u$$

admits homogenization under the natural condition $\langle c \rangle = 0$ ($\langle \cdot \rangle$ stands for the mean value) and that the homogenized operator takes the form

$$\frac{\partial}{\partial t} u = \frac{1}{2} \hat{a} \frac{\partial^2}{\partial x^2} u + \hat{c} u,$$

with constant $\hat{a}$ and $\hat{c}$.

In contrast with symmetric divergence form parabolic problems, in the presence of the lower order terms the asymptotic behaviour of operators with random coefficients might differ a lot from that of periodic operators.

Homogenization problem for parabolic operators whose coefficients are periodic in spatial variables and random stationary in time, were studied in [4, 5, 14]. It was shown that, under natural mixing assumptions on the coefficients, the critical rate of the potential growth is of order $1/\varepsilon$. In this case the limit equation is a stochastic PDE.

If the oscillating potential is random stationary (statistically homogeneous) in spatial variables then the range of the oscillations (the power of $\varepsilon^{-1}$ in front of potential $c$) should depend on the spatial dimension.

In this work we deal with a one-dimensional spatial variable and show that the range of oscillation should be of order $\frac{1}{\sqrt{\varepsilon}}$. This means that for larger powers of $\frac{1}{\varepsilon}$ the family of solutions is not tight as $\varepsilon \to 0$, while for smaller powers of $\frac{1}{\varepsilon}$ the contribution of the potential is asymptotically negligible.

It turns out that the Dirichlet forms technique which is usually quite efficient in homogenization problems, does not apply to problem (1.1) because one cannot prove any lower bound for the quadratic form corresponding to the operator (1.1). This is due to the fact that the problem is stated on the whole line $\mathbb{R}$, and not on a compact interval, and that the coefficients of the operator are a.s. unbounded, see the discussion in Section 6. Instead we use the direct approach combining the Feynman–Kac formula with several correctors, Itô calculus and martingale convergence arguments.

The main result of the paper (see Theorem 2.2) states that under proper mixing conditions the solution $u^\varepsilon$ of eq. (1.1) converges in law to a random field

$$u(x,t) = \mathbb{E}\bigg( g(x + \sqrt{\tilde{a}} B_t) \exp\bigg( \frac{\overline{c}}{\tilde{a}} \int_\mathbb{R} L_t^{y-x} W(\mathrm{d}y) \bigg) \bigg),$$

where $B_\cdot$ and $W_\cdot$ are independent Brownian motions, $\mathbb{E}$ and $L_t^y$ are respectively the expectation and the local time related to $\sqrt{\tilde{a}} B_\cdot$, and $\tilde{a}$, $\overline{c}$ are constants.

The interpretation of this expression is given in the last section of the paper. It is shown that the effective equation is not a standard SPDE but rather a parabolic PDE with random coefficients.

Let us give an intuitive explanation of our result. The Feynman–Kac formula for the solution of eq. (1.1) yields

$$u^\varepsilon(t,x) = \mathbb{E}\bigg[ g(X_t^{\varepsilon,x}) \exp\bigg( \frac{1}{\sqrt{\varepsilon}} \int_0^t c\bigg(\frac{X_s^{\varepsilon,x}}{\varepsilon}\bigg) \mathrm{d}s \bigg) \bigg],$$

where $\mathbb{E}$ means expectation with respect to the law of the diffusion $X_\cdot^{\varepsilon,x}$, the random field $c$ being frozen (or "quenched"). Under the assumptions which we shall make below, one can apply a version of the functional central limit theorem, which tells us that

$$W_\varepsilon(x) = \frac{1}{\sqrt{\varepsilon}} \int_0^x c\bigg(\frac{y}{\varepsilon}\bigg) \mathrm{d}y$$

converges weakly towards $\overline{c} W(x)$, where $W$ is a standard Wiener process. Now the exponent in the above Feynman–Kac formula reads

$$\int_0^t W_\varepsilon'(X_s^{\varepsilon,x}) \mathrm{d}s = \int_\mathbb{R} W_\varepsilon'(y) L_t^{y,\varepsilon} \mathrm{d}y,$$

where $L_t^{y,\varepsilon}$ denotes the local time at time $t$ and point $x$ of the diffusion process $\{X^{\varepsilon,x}\}$. One might expect that the last integral converges towards the integral of the limiting local time, with respect to the limiting Wiener process. This is one of the results which will be established in this paper.

Our paper is organized as follows. In Section 2, we formulate our assumptions, and the results. In Section 3, we prove some weak convergence results. Section 4 is devoted to the proof of the pointwise convergence of the sequence $u^\varepsilon(t,x)$, while Section 5 is concerned with convergence in the space of continuous functions. Finally in Section 6 we discuss the limiting PDE.

## 2. Set up and statement of the main result

We make the following assumptions:

(A.1) The initial condition $g$ belongs to $L^2(\mathbb{R}) \cap \mathcal{C}_b(\mathbb{R})$.

(A.2) The coefficients $\{a(x), x \in \mathbb{R}\}$ and $\{c(x), x \in \mathbb{R}\}$ are stationary random fields defined on a probability space $(\Sigma, \mathcal{A}, P)$, and we assume that

$$0 < c \le a(x) \le C, \quad x \in \mathbb{R}, P \text{ a.s.}, \tag{2.1}$$

$$Ec(0) = 0, \quad \|c(0)\|_{L^\infty(\Sigma)} < \infty, \tag{2.2}$$

$$\int_{-\infty}^{\infty} |Ec(0)c(x)| \, dx < \infty, \tag{2.3}$$

where $E$ denotes expectation with respect to the probability measure $P$.

(A.3) Let

$$\mathcal{F}_x := \sigma\{a(y), c(y); y \le x\}, \qquad \mathcal{F}^x := \sigma\{a(y), c(y); y \ge x\}.$$

We assume that the random fields $a$ and $c$ are $\phi$-mixing in the following sense. Define, for $h > 0$, $\phi(h)$ the mixing coefficient with respect to the $\sigma$-algebras from above, as

$$\phi(h) := \sup_{\{A \in \mathcal{F}_x, B \in \mathcal{F}^{x+h}, P(A) > 0\}} |P(B|A) - P(B)|.$$

We suppose that

$$\int_0^\infty \phi^{1/2}(h) \, dh < \infty. \tag{2.4}$$

Consider now the family of Dirichlet forms $\{\mathcal{E}^{\varepsilon,\sigma}, \varepsilon > 0, \sigma \in \Sigma\}$ on $L^2(\mathbb{R})$ defined by

$$\mathcal{E}^{\varepsilon,\sigma}(u,v) = \frac{1}{2} \int_\mathbb{R} a\left(\frac{x}{\varepsilon}, \sigma\right) u'(x) v'(x) \, dx,$$

with domain $H^1(\mathbb{R})$. For each $\varepsilon > 0, \sigma \in \Sigma$ there exists a unique self-adjoint operator $L^{\varepsilon,\sigma}$ with domain $\mathcal{D}(L^{\varepsilon,\sigma})$, such that

$$\mathcal{E}_\varepsilon^\sigma(u,v) = -(L^{\varepsilon,\sigma} u, v)_{L^2(\mathbb{R})},$$

for $u \in \mathcal{D}(L^{\varepsilon,\sigma}), v \in L^2(\mathbb{R})$.

For each initial point $x \in \mathbb{R}$, $\sigma \in \Sigma$ and $\varepsilon > 0$, there exists a continuous Markov process $\{X_t^{\varepsilon,\sigma,x}, t \ge 0\}$ defined on some probability space $(\Omega, \mathcal{F}, \mathbb{P}^{\varepsilon,x,\sigma})$ whose generator is $L^{\varepsilon,\sigma}$ and which starts at time $t = 0$ from $x$. The probability may depend on the three parameters $\varepsilon$, $x$, $\sigma$. $x$ will be fixed throughout this paper, so we drop it from now on. Note that the process $\{X_t^{\varepsilon,x}, t \ge 0\}$ is in fact defined on the probability space $(\Sigma \times \Omega, \mathcal{A} \otimes \mathcal{F}, Q^\varepsilon)$ (as such, it is not a Markov process), where the probability $Q^\varepsilon$ on the product space $\Sigma \times \Omega$ is defined as

$$Q^\varepsilon(A) = \int \int_A P(\mathrm{d}\sigma) \mathbb{P}^{\varepsilon,\sigma}(\mathrm{d}\omega).$$

The Feynman–Kac formula allows us to write down an explicit formula for the solution of eq. (1.1):

$$u^\varepsilon(t,x) = \mathbb{E}^{\varepsilon,\cdot}\left[g(X_t^{\varepsilon,x}) \exp\left(\frac{1}{\sqrt{\varepsilon}} \int_0^t c\left(\frac{X_s^{\varepsilon,x}}{\varepsilon}\right) \mathrm{d}s\right)\right], \tag{2.5}$$

where $\mathbb{E}^{\varepsilon,\sigma}$ denotes expectation with respect to $\mathbb{P}^{\varepsilon,\sigma}$.

Considering assumption (A.2), we define the finite quantities

$$\overline{c}^2 = \int_{-\infty}^{\infty} E[c(0)c(x)]\, \mathrm{d}x, \qquad \tilde{a} = \left[E\left(\frac{1}{a(x)}\right)\right]^{-1}. \tag{2.6}$$

In view of Theorem 5.1 and Lemma 5.1 from [11], we may state the following theorem.

**Theorem 2.1.** *We have the following convergence, P a.s.:*

$$X_\cdot^{\varepsilon,x} \Rightarrow X_\cdot^x := x + X_\cdot,$$

*in $\mathcal{C}([0,\infty))$, where $X$ is a one-dimensional Brownian motion defined on $(\Omega, \mathcal{F}, \mathbb{P})$ such that $\mathbb{E}(X_t^2) = \tilde{a}t$, for any $t \geq 0$.*

The main result of this paper is the following theorem.

**Theorem 2.2.** *Let*

$$u(t,x) := \mathbb{E}\left[g(x + X_t) \exp\left(\frac{\overline{c}}{\tilde{a}} \int_\mathbb{R} L_t^{y-x} W(\mathrm{d}y)\right)\right],$$

*where $W$ denotes a one-dimensional standard Brownian motion defined on the probability space $(\Sigma, \mathcal{A}, P)$ and $L_t^y$ is the local time at time $t$ and point $y$ of the process $\{X_t, t \geq 0\}$ defined on $(\Omega, \mathcal{F}, \mathbb{P})$.*
*Then $u^\varepsilon \Rightarrow u$ in law in $C(\mathbb{R}_+ \times \mathbb{R})$, as $\varepsilon \to 0$.*

We introduce the notation

$$Y_t^{\varepsilon,x} := \frac{1}{\sqrt{\varepsilon}} \int_0^t c\left(\frac{X_s^{\varepsilon,x}}{\varepsilon}\right) \mathrm{d}s.$$

The first step in the proof of Theorem 2.2 is to establish the weak convergence of the pair $(X_t^{\varepsilon,x}, Y_t^{\varepsilon,x})$, which is done in the next section.

## 3. Weak convergence

The main result of this section is the following theorem.

**Theorem 3.1.** *For each $t > 0$,*

$$(X_t^{\varepsilon,x}, Y_t^{\varepsilon,x}) \Rightarrow (X_t^x, Y_t^x)$$

*weakly, as $\varepsilon \to 0$, with*

$$Y_t^x := \frac{\overline{c}}{\tilde{a}} \int_{\mathbb{R}} L_t^{y-x} W(\mathrm{d}y),$$

*where, as above, $L_t^y$ is the local time at point $y$ and time $t$ of the Brownian motion $\{X_t, t \geq 0\}$ defined on $(\Omega, \mathcal{F}, \mathbb{P})$, and $\{W_y, y \in \mathbb{R}\}$ is a Wiener process defined on $(\Sigma, \mathcal{A}, P)$, so that $(X, L)$ and $W$ are independent.*

Theorem 3.1 will follow easily from Propositions 3.7 and 3.10, as we shall see at the end of this section. Note that all we shall need in the next section is both Propositions 3.7 and 3.10.

Let us first state a consequence of Aronson's estimate, see Lemma II.1.2 in [16]:

**Lemma 3.2.** *There exists $\kappa > 0$, which depends only on $c$ and $C$ in (2.1), such that for all $\varepsilon > 0$, $r > 0$,*

$$\mathbb{P}\left(\sup_{0 \leq s \leq t} |X_s^{\varepsilon,x} - x| > r\right) \leq \kappa \exp\left(\frac{-r^2}{\kappa t}\right).$$

We next prove the easiest part of the above result, i.e. we give a proof of Theorem 2.1, since we shall need some of its details later.

**Proof of Theorem 2.1.** Let $\{\chi(x), x \in \mathbb{R}\}$ be the zero mean random process given by the formula

$$\chi(x) = \tilde{a} \int_0^x \frac{\mathrm{d}y}{a(y)} - x.$$

We note that from Birkhof's ergodic theorem (see e.g. Theorem 24.1 in [3]),

$$\frac{\chi(x)}{x} \to 0, \quad P \text{ a.s., as } |x| \to \infty. \tag{3.1}$$

Moreover, this random process satisfies the two relations:

$$(a(1+\chi'))'(x) = 0, \quad x \in \mathbb{R}, \tag{3.2}$$

and

$$a(x)(1+\chi'(x)) = \tilde{a}, \quad x \in \mathbb{R}.$$

We now define

$$Z_t^\varepsilon = X_t^{\varepsilon,x} + \varepsilon \chi\left(\frac{X_t^{\varepsilon,x}}{\varepsilon}\right).$$

It follows from the Itô–Fukushima decomposition (see [7] or Theorem 0.10 in [11]) and (3.2) that (here and further below, $M^{X^{\varepsilon,x}}$ denotes the martingale part of the process $X^{\varepsilon,x}$)

$$Z_t^\varepsilon = Z_0^\varepsilon + \int_0^t \left[1 + \chi'\left(\frac{X_s^{\varepsilon,x}}{\varepsilon}\right)\right] \mathrm{d}M_s^{X^{\varepsilon,x}}, \quad t \geq 0,$$

hence $\{Z^\varepsilon\}$ is $P$ a.s. a $\mathbb{P}$-martingale. Moreover, its quadratic variation is given by

$$\langle Z^\varepsilon \rangle_t = \tilde{a}^2 \int_0^t \frac{\mathrm{d}s}{a(X_s^{\varepsilon,x}/\varepsilon)}.$$

It will be proved below in Lemma 3.9 that

$$\langle Z^\varepsilon \rangle_t \to \tilde{a}t$$

in $Q^\varepsilon$ probability. It now follows from well-known results that $P$ a.s.,

$$Z^\varepsilon \Rightarrow x + \sqrt{\tilde{a}}B,$$

where $\{B_t, t \geq 0\}$ is a standard Brownian motion defined on the probability space $(\Omega, \mathcal{F}, \mathbb{P})$.

Moreover, for all $T > 0$,

$$\sup_{0 \leq t \leq T} |X_t^{\varepsilon,x} - Z_t^\varepsilon| \to 0, \quad Q^\varepsilon \text{ a.s.,}$$

consequently $P$ a.s.,

$$X^{\varepsilon,x} \Rightarrow x + \sqrt{\tilde{a}}B$$

in $\mathbb{P}$ law. □

Let $\Phi$ denote the solution of the ordinary differential equation:

$$(a\Phi')'(x) = c(x),$$

which is defined as follows:

$$\Phi'(x) = \frac{1}{a(x)} \int_0^x c(y) \, dy,$$

$$\Phi(x) = \int_0^x \left( \frac{1}{a(z)} \int_0^z c(y) \, dy \right) dz. \tag{3.3}$$

We let

$$\bar{c}W_\varepsilon(x) := \sqrt{\varepsilon} a\left(\frac{x}{\varepsilon}\right) \Phi'\left(\frac{x}{\varepsilon}\right)$$

$$= \frac{1}{\sqrt{\varepsilon}} \int_0^x c\left(\frac{y}{\varepsilon}\right) dy \tag{3.4}$$

and

$$F_\varepsilon(x) := \varepsilon^{3/2} \Phi\left(\frac{x}{\varepsilon}\right)$$

$$= \frac{1}{\sqrt{\varepsilon}} \int_0^x \frac{1}{a(z/\varepsilon)} \int_0^z c\left(\frac{y}{\varepsilon}\right) dy \, dz$$

$$= \bar{c} \int_0^x \frac{1}{a(z/\varepsilon)} W_\varepsilon(z) \, dz. \tag{3.5}$$

We first prove the following proposition.

**Proposition 3.3.** *The sequence of random processes $\{W_\varepsilon\}$ converges weakly, as $\varepsilon \to 0$, in the space $\mathcal{C}(\mathbb{R})$, to a standard Wiener process $\{W\}$ defined on $(\Sigma, \mathcal{A}, P)$.*

**Proof.** Denote, for $x \geq 0$, $W_\varepsilon^1(x) = W_\varepsilon(x)$ and $W_\varepsilon^2(x) = W_\varepsilon(-x)$. According to assumptions (A.1), (A.2), (A.3) and the functional central limit theorem (see e.g. [2], pages 178, 179), it follows that

$$(W_\varepsilon^1, W_\varepsilon^2) \xrightarrow{\mathcal{D}} (W^1, W^2),$$

where $\{W^1(x), x \geq 0\}$ and $\{W^2(x), x \geq 0\}$ are mutually independent standard Brownian motions. Finally we denote by $\{W(x), x \in \mathbb{R}\}$ the process defined by

$$W(x) := W^1(x), \quad \text{for } x \geq 0, \qquad W(x) := W^2(-x), \quad \text{for } x < 0. \qquad \square$$

It remains to show why Theorem 3.1 follows from Theorem 2.1 and Proposition 3.3.
First we define, for $x \in \mathbb{R}$,

$$k_\varepsilon(x) = \int_0^x \frac{1}{a(z/\varepsilon)} \, \mathrm{d}z, \qquad k(x) = \frac{x}{\tilde{a}}.$$

We have the following lemma.

**Lemma 3.4.** *As $\varepsilon \to 0$,*

$$k_\varepsilon \to k \quad \text{in } C(\mathbb{R}), P \text{ a.s.}$$

**Proof.** Since $a(\cdot)$ is bounded away from zero, $a^{-1}(\cdot)$ is bounded. Hence the collection of random functions $k_\varepsilon$ is tight in $C(\mathbb{R})$. It then suffices to show that the finite dimensional marginals converge in law to those of the deterministic function $k$. But from Birkhoff's ergodic theorem, for any $x_1, x_2, \ldots, x_n \in \mathbb{R}$,

$$(k_\varepsilon(x_1), k_\varepsilon(x_2), \ldots, k_\varepsilon(x_n)) \to \left(\frac{x_1}{\tilde{a}}, \frac{x_2}{\tilde{a}}, \ldots, \frac{x_n}{\tilde{a}}\right)$$

in $P$ a.s., as $\varepsilon \to 0$. $\qquad \square$

Denote by $C_+(\mathbb{R})$ the space of continuous and increasing functions on $\mathbb{R}$, and by $S$ the map from $\mathcal{C}(\mathbb{R}) \times C_+(\mathbb{R})$ into $\mathcal{C}(\mathbb{R})$ defined by

$$S((f,h))(x) := \int_0^x f(z) \, \mathrm{d}h(z).$$

We have the following lemma.

**Lemma 3.5.** *The mapping $S$ is continuous, from $E = C(\mathbb{R}) \times C_+(\mathbb{R})$, equipped with the product of the locally uniform topology of $C(\mathbb{R}) \times C(\mathbb{R})$, into $C(\mathbb{R})$, equipped with the locally uniform topology.*

**Proof.** It suffices to show that for each $N > 0$, if $\{(f_n, h_n)\} \subset C([-N, N]) \times C_+([-N, N])$, and

$$\sup_{|x| \leq N} (|f_n(x) - f(x)| + |h_n(x) - h(x)|) \to 0, \quad \text{as } n \to \infty,$$

then

$$\sup_{|x| \leq N} |S(f_n, h_n)(x) - S(f, h)(x)| \to 0, \quad \text{as } n \to \infty.$$

But this follows from Lemma 5.8 in [8]. $\qquad \square$

We now have the following lemma.

**Lemma 3.6.** *As $\varepsilon \to 0$,*

$$(W_\varepsilon, F_\varepsilon) \to (W, F)$$

*in $C(\mathbb{R}) \times C(\mathbb{R})$ in law, where $F(x) = \frac{\bar{c}}{\tilde{a}} \int_0^x W(z) \, \mathrm{d}z$ and $W_\varepsilon$ and $F_\varepsilon$ are defined in (3.4) and (3.5) respectively.*

**Proof.** Since $k_\varepsilon$ converges to a deterministic limit, it follows from a well-known theorem (see e.g. Theorem 4.4 in [2]) that the pair $(W_\varepsilon, k_\varepsilon)$ converges. Hence from Lemma 3.5, the pair $(W_\varepsilon, F_\varepsilon) = (W_\varepsilon, S(W_\varepsilon, k_\varepsilon))$ converges. □

The next step is to show that the triple $(X_\varepsilon, W_\varepsilon, F_\varepsilon)$ converges. This is essentially a consequence of the three following facts: $X_\varepsilon$ converges, $(W_\varepsilon, F_\varepsilon)$ converges, and the two limits $X$ and $(W, F)$ are defined on $(\Omega, \mathcal{F}, \mathbb{P})$ and $(\Sigma, \mathcal{A}, P)$ respectively. We now prove that fact rigorously.

**Proposition 3.7.** *For any $x \in \mathbb{R}$, as $\varepsilon \to 0$,*

$$(X^{\varepsilon,x}, W_\varepsilon, F_\varepsilon) \to (X^x, W, F)$$

*in law, in $C(\mathbb{R}_+) \times C(\mathbb{R}) \times C(\mathbb{R})$.*

**Proof:.** We first choose two arbitrary functionals $\Psi \in \mathcal{C}_b(\mathcal{C}(\mathbb{R}) \times \mathcal{C}(\mathbb{R}))$ and $\Theta \in \mathcal{C}_b(\mathcal{C}(\mathbb{R}_+))$. We have

$$\int_{\Sigma \times \Omega} \Psi(W_\varepsilon(\sigma), F_\varepsilon(\sigma)) \Theta(X^{\varepsilon,x}(\omega, \sigma)) P(d\sigma) \mathbb{P}^{\varepsilon,\sigma}(d\omega)$$

$$= \int_\Sigma \Psi(W_\varepsilon(\sigma), F_\varepsilon(\sigma)) \left( \int_\Omega \Theta(X^{\varepsilon,x}(\omega, \sigma)) \mathbb{P}^{\varepsilon,\sigma}(d\omega) \right) P(d\sigma)$$

$$= \int_\Sigma \Psi(W_\varepsilon(\sigma), F_\varepsilon(\sigma)) \left[ \int_\Omega (\Theta(X^{\varepsilon,x}(\omega, \sigma)) - \Theta(X^x(\omega))) \mathbb{P}^{\varepsilon,\sigma}(d\omega) \right] P(d\sigma)$$

$$+ \int_\Omega \Theta(X^x(\omega)) \mathbb{P}^{\varepsilon,\sigma}(d\omega) \times \int_\Sigma \Psi(W_\varepsilon(\sigma), F_\varepsilon(\sigma)) P(d\sigma).$$

Theorem 2.1 tells us that

$$\int_\Omega \Theta(X^{\varepsilon,x}(\omega, \sigma)) \mathbb{P}^{\varepsilon,\sigma}(d\omega) \to \int_\Omega \Theta(X^x(\omega)) \mathbb{P}(d\omega),$$

$P$ a.s., and from Lemma 3.6,

$$\int_\Sigma \Psi(W_\varepsilon(\sigma), F_\varepsilon(\sigma)) P(d\sigma) \to \int_\Sigma \Psi(W(\sigma), F(\sigma)) P(d\sigma),$$

as $\varepsilon \to 0$. Hence, from the Bounded Convergence Theorem, we conclude that

$$\int_{\Sigma \times \Omega} \Psi(W_\varepsilon(\sigma), F_\varepsilon(\sigma)) \Theta(X^{\varepsilon,x}(\omega, \sigma)) P(d\sigma) \mathbb{P}^{\varepsilon,\sigma}(d\omega)$$

$$\to \int_\Omega \Theta(X^x(\omega)) P(d\omega) \times \int_\Sigma \Psi(W(\sigma), F(\sigma)) P(d\sigma).$$

It now suffices to note that $\mathcal{A} := \{\Psi \otimes \Theta, \Psi \in \mathcal{C}_b(\mathcal{C}(\mathbb{R}) \times \mathcal{C}(\mathbb{R})), \Theta \in \mathcal{C}_b(\mathcal{C}(\mathbb{R}_+))\}$ is a determining class on $\mathcal{C}(\mathbb{R}) \times \mathcal{C}(\mathbb{R}) \times \mathcal{C}(\mathbb{R}_+)$. □

It follows from Lemma 3.8 that

$$\varepsilon^{3/2} \Phi\left(\frac{X_t^{\varepsilon,x}}{\varepsilon}\right) = \varepsilon^{3/2} \Phi\left(\frac{x}{\varepsilon}\right) + \frac{1}{2\sqrt{\varepsilon}} \int_0^t c\left(\frac{X_s^{\varepsilon,x}}{\varepsilon}\right) ds + M_t^{\varepsilon,x}, \tag{3.6}$$

in other words

$$F_\varepsilon(X_t^{\varepsilon,x}) = F_\varepsilon(x) + \frac{1}{2} Y_t^{\varepsilon,x} + M_t^{\varepsilon,x},$$

where $\{M_t^{\varepsilon,x}, t \geq 0\}$ is the continuous martingale

$$M_t^{\varepsilon,x} = \sqrt{\varepsilon} \int_0^t \Phi'\left(\frac{X_s^{\varepsilon,x}}{\varepsilon}\right) dM_s^{X^{\varepsilon,x}},$$

and $M^{X^{\varepsilon,x}}$ denotes again the martingale part of the process $X^{\varepsilon,x}$. In particular the quadratic variation of $M^{\varepsilon,x}$ is given by the quantity

$$\langle M^{\varepsilon,x} \rangle_t = \varepsilon \int_0^t a\left(\frac{X_s^{\varepsilon,x}}{\varepsilon}\right) \left[\Phi'\left(\frac{X_s^{\varepsilon,x}}{\varepsilon}\right)\right]^2 ds$$

$$= \int_0^t \frac{1}{a(X_s^{\varepsilon,x}/\varepsilon)} \left[\sqrt{\varepsilon}\Phi'\left(\frac{X_s^{\varepsilon,x}}{\varepsilon}\right) a\left(\frac{X_s^{\varepsilon,x}}{\varepsilon}\right)\right]^2 ds$$

$$= \bar{c}^2 \int_0^t \frac{1}{a(X_s^{\varepsilon,x}/\varepsilon)} (W_\varepsilon(X_s^{\varepsilon,x}))^2 ds,$$

and the joint quadratic variation of $M^{\varepsilon,x}$ and $Z^\varepsilon$ is

$$\langle M^{\varepsilon,x}, Z^\varepsilon \rangle_t = \sqrt{\varepsilon} \int_0^t \Phi'\left(\frac{X_s^{\varepsilon,x}}{\varepsilon}\right) \left[1+\chi'\left(\frac{X_s^{\varepsilon,x}}{\varepsilon}\right)\right] a\left(\frac{X_s^{\varepsilon,x}}{\varepsilon}\right) ds$$

$$= \tilde{a}\bar{c} \int_0^t \frac{W_\varepsilon(X_s^{\varepsilon,x})}{a(X_s^{\varepsilon,x}/\varepsilon)} ds.$$

**Lemma 3.8.** *The identity* (3.6) *holds $Q^\varepsilon$ a.s.*

**Proof.** Let $\{X_t^{\varepsilon,x,M}, t \geq 0\}$ denote the process $\{X_t^{\varepsilon,x}, t \geq 0\}$, killed when exiting the interval $[-\varepsilon(M+1), \varepsilon(M+1)]$, and $L^{\varepsilon,\sigma,M}$ its generator. For any $M > 0$, the random function

$$\Phi_M(x) = \begin{cases} \int_0^x \frac{1}{a(z)} \int_0^z c_M(y)\,dy\,dz, & \text{if } x \geq 0; \\ \int_x^0 \frac{1}{a(z)} \int_z^0 c_M(y)\,dy\,dz, & \text{if } x < 0, \end{cases}$$

with the random field $\{c_M(x), -M-1 \leq x \leq M+1\}$ defined by

$$c_M(x) = \begin{cases} \beta_M, & \text{if } -M-1 < x < -M; \\ c(x), & \text{if } -M \leq x \leq M; \\ \alpha_M, & \text{if } M < x < M+1, \end{cases}$$

where

$$\alpha_M = -\frac{\int_0^{M+1}(1/a(z)) \int_0^{z\wedge M} c(y)\,dy\,dz}{\int_M^{M+1}(z-M)/a(z)\,dz};$$

$$\beta_M = -\frac{\int_{-M-1}^0 (1/a(z)) \int_{z\vee -M}^0 c(y)\,dy\,dz}{\int_{-M-1}^{-M}(z+M)/a(z)\,dz}$$

satisfies

$$\begin{cases} (L^{\varepsilon,\sigma}\Phi_M)(\frac{x}{\varepsilon}) = \varepsilon^{-2} c_M(\frac{x}{\varepsilon}), \\ \Phi_M(-M-1) = \Phi_M(M+1) = 0. \end{cases}$$

Consequently for each $\varepsilon > 0$,

$$\Phi_M\left(\frac{\cdot}{\varepsilon}\right) \in \mathcal{D}(L^{\varepsilon,\cdot,M}),$$

and by virtue of the Itô–Fukushima decomposition (see [7] or Theorem 0.10 in [11]), we get for $|x| \leq M+1$

$$\varepsilon^{3/2}\Phi_M\left(\frac{X^{\varepsilon,x}_{t\wedge\tau^{\varepsilon,x}_M}}{\varepsilon}\right) = \varepsilon^{3/2}\Phi_M\left(\frac{x}{\varepsilon}\right) + \frac{1}{2\sqrt{\varepsilon}}\int_0^{t\wedge\tau^{\varepsilon,x}_M} c_M\left(\frac{X^{\varepsilon,x}_s}{\varepsilon}\right)\mathrm{d}s + M^{\varepsilon,x}_{t\wedge\tau^{\varepsilon,x}_M},$$

where

$$\tau^{\varepsilon,x}_M = \inf\{t \geq 0, |X^{\varepsilon,x}_t| = \varepsilon(M+1)\}.$$

The result follows, by letting $M \to \infty$, with the help of Lemma 3.2. $\square$

Define

$$h_\varepsilon(t) := \int_0^t \frac{1}{a(X^{\varepsilon,x}_s/\varepsilon)}\,\mathrm{d}s.$$

**Lemma 3.9.** *As $\varepsilon \to 0$,*

$$h_\varepsilon(\cdot) \to \frac{\cdot}{\tilde{a}}$$

*in $\mathcal{C}(\mathbb{R}_+)$ in probability.*

**Proof.** Denote $\theta(x) = \frac{1}{a(x)} - \frac{1}{\tilde{a}}$. Then $\theta(x)$ is a bounded stationary field with zero mean. Letting

$$\Theta(x) = \int_0^x \left(\frac{1}{a(z)}\int_0^z \theta(y)\,\mathrm{d}y\right)$$

and repeating the argument in Lemma 3.8, we get

$$\varepsilon^{3/2}\Theta\left(\frac{X^{\varepsilon,x}_t}{\varepsilon}\right) = \varepsilon^{3/2}\Theta\left(\frac{x}{\varepsilon}\right) + \frac{1}{\sqrt{\varepsilon}}\int_0^t \theta\left(\frac{X^{\varepsilon,x}_s}{\varepsilon}\right)\mathrm{d}s + \mathcal{M}^\varepsilon_t, \tag{3.7}$$

where

$$\langle \mathcal{M}^\varepsilon \rangle_t = \varepsilon \int_0^t a^{-1}\left(\frac{X^{\varepsilon,x}_s}{\varepsilon}\right)\left(\int_0^{X^{\varepsilon,x}_s/\varepsilon}\theta(y)\,\mathrm{d}y\right)^2 \mathrm{d}s.$$

In the same way as in the proof of Proposition 3.7, one can show that the families $\{\varepsilon^{3/2}\Theta(\frac{X^{\varepsilon,x}_\cdot}{\varepsilon}) - \varepsilon^{3/2}\Theta(\frac{x}{\varepsilon})\}$ and $\{\mathcal{M}^\varepsilon_\cdot\}$ are tight in $\mathcal{C}(\mathbb{R})$. Indeed, $\Theta$ is constructed from $\theta$ exactly as $\Phi$ from $c$. Moreover, $\theta$ is, exactly as $c$, a stationary mixing bounded and zero mean random field. Now, multiplying the relation (3.7) by $\sqrt{\varepsilon}$, we conclude that

$$\int_0^\cdot \theta\left(\frac{X^{\varepsilon,x}_s}{\varepsilon}\right)\mathrm{d}s \to 0$$

in $\mathcal{C}(\mathbb{R}_+)$ in probability. This implies the desired convergence. $\square$

For any $f \in C_b(\mathbb{R})$ and $\varepsilon > 0$, let us define the process $\{N_t^{f,\varepsilon}; t \geq 0\}$ by

$$N_t^{f,\varepsilon} = \int_0^t \frac{f(X_s^{\varepsilon,x})}{a(X_s^{\varepsilon,x}/\varepsilon)} \, dM_s^{X^{\varepsilon,x}}$$

$$= \frac{1}{\tilde{a}} \int_0^t f(X_s^{\varepsilon,x}) \, dZ_s^{\varepsilon}.$$

We now prove the following proposition.

**Proposition 3.10.** *$P$ a.s.,*

$$(N^{f,\varepsilon}, X^{\varepsilon,x}) \Rightarrow (N^f, X^x) \quad \text{in } C([0,t]) \times C([0,t]),$$

*where*

$$N_s^f = \frac{1}{\tilde{a}} \int_0^s f(X_r^x) \, dX_r^x.$$

**Proof.** Since $(X^{\varepsilon,x} - Z^{\varepsilon})$ converges to zero in probability, $(N_t^{f,\varepsilon}, X_t^{\varepsilon,x})$ behaves as $\varepsilon \to 0$ exactly as $(N_t^{f,\varepsilon}, Z_t^{\varepsilon})$, hence we consider the two-dimensional martingale $(N_t^{f,\varepsilon}, Z_t^{\varepsilon})$, and compute its associated bracket process, which takes values in the set of $2 \times 2$ symmetric matrices. We have

$$\left\langle\!\!\left\langle \begin{pmatrix} N^{f,\varepsilon} \\ Z^{\varepsilon} \end{pmatrix} \right\rangle\!\!\right\rangle_t = \begin{pmatrix} \int_0^t \frac{f^2(X_s^{\varepsilon,x})}{a(X_s^{\varepsilon,x}/\varepsilon)} \, ds & \tilde{a} \int_0^t \frac{f(X_s^{\varepsilon,x})}{a(X_s^{\varepsilon,x}/\varepsilon)} \, ds \\ \tilde{a} \int_0^t \frac{f(X_s^{\varepsilon,x})}{a(X_s^{\varepsilon,x}/\varepsilon)} \, ds & \tilde{a}^2 \int_0^t \frac{ds}{a(X_s^{\varepsilon,x}/\varepsilon)} \end{pmatrix}.$$

Combining Theorem 2.1, Lemmas 3.5 and 3.9 we obtain that this $\mathbb{R}^4$-valued process converges $P$ a.s. in $\mathbb{P}$ law towards

$$\begin{pmatrix} \frac{1}{\tilde{a}} \int_0^t f^2(X_s^x) \, ds & \int_0^t f(X_s^x) \, ds \\ \int_0^t f(X_s^x) \, ds & \tilde{a} t \end{pmatrix}.$$

We then conclude that $P$ a.s.,

$$(N^{f,\varepsilon}, X^{\varepsilon,x}) \Rightarrow (N^f, X^x)$$

in $\mathbb{P}$ law, where

$$N_t^f = \frac{1}{\tilde{a}} \int_0^t f(X_s^x) \, dX_s^x, \quad t \geq 0. \qquad \square$$

The statement below is a straightforward consequence of Birkhoff's ergodic theorem.

**Proposition 3.11.** *For any $N > 0$, $x \in \mathbb{R}$ and $f \in C(\mathbb{R})$ the following convergence holds $P$ a.s.*

$$\sup_{|y-x| \leq N} \left| \int_x^y \frac{f(z)}{a(z/\varepsilon)} \, dz - \frac{1}{\tilde{a}} \int_x^y f(z) \, dz \right| \longrightarrow 0.$$

We now establish the version of (3.6) for $\varepsilon = 0$, which is an Itô type formula for the process $\{F(x + X_t), t \geq 0\}$, where $F(y) := \frac{\overline{c}}{\tilde{a}} \int_0^y W(z) \, dz, y \in \mathbb{R}$. More precisely, we have the following lemma.

**Lemma 3.12.** *For any $t \geq 0$ and $x \in \mathbb{R}$, we have*

$$F(X_t^x) - F(x) = \frac{\overline{c}}{\tilde{a}} \int_0^t W(X_s^x) \, dX_s + \frac{\overline{c}}{2} \int_{\mathbb{R}} L_t^{y-x} W(dy). \tag{3.8}$$

**Proof.** We prove this formula by using smooth approximations of the process $\{W\}$, obtained by convolution. Let $\rho$ be a $\mathcal{C}_0^2(\mathbb{R})$ function such that $\rho \geq 0$, $\text{supp}(\rho) \subseteq [-1, 1]$ and $\int_\mathbb{R} \rho(y)\,dy = 1$. Define now

$$\rho_n(y) := n\rho(ny), \quad W_n(y) := \int_\mathbb{R} \rho_n(y-z) W(z)\,dz = \int_{-1}^1 \rho(z) W\left(y - \frac{z}{n}\right) dz. \tag{3.9}$$

From the uniform continuity of $W$ on compacts, $\|W_n - W\|_{\mathcal{C}(K)} \to 0$, $P$ a.s., as $n \to \infty$, for any compact set $K$ in $\mathbb{R}$. Moreover, taking into account the fact that $W$ is a standard Brownian motion, we get

$$E(W_n^4(y)) \leq C \int_{-1}^1 \left[ E\left(W^4\left(y - \frac{z}{n}\right)\right) \right] dz = C \int_{-1}^1 \left(y - \frac{z}{n}\right)^2 dz \leq C(1+y^2), \tag{3.10}$$

for $y \in \mathbb{R}$. Set

$$F_n(\cdot) := \frac{\overline{c}}{\tilde{a}} \int_0^\cdot W_n(y)\,dy.$$

Itô's formula applied to the process $\{F_n(x + X_t), t \geq 0\}$ gives

$$F_n(x + X_t) - F_n(x) = \frac{\overline{c}}{\tilde{a}} \int_0^t W_n(x + X_s)\,dX_s + \frac{\overline{c}}{2} \int_0^t W_n'(x + X_s)\,ds. \tag{3.11}$$

Recall that $\{x + X_t, t \geq 0\}$ is a non-standard Brownian motion independent of $\{W\}$. It is easy to see that the left-hand side in the last formula tends to $F(x + X_t) - F(x)$, $P \times \mathbb{P}$ a.s., as $n \to \infty$. Moreover,

$$E \times \mathbb{E}\left\{\left[\int_0^t (W_n(x+X_s) - W(x+X_s))\,dX_s\right]^2\right\}$$

$$= \tilde{a} E\left\{\mathbb{E}\left[\int_0^t (W_n(x+X_s) - W(x+X_s))^2\,ds\right]\right\}$$

$$\leq \tilde{a} \mathbb{E}\left\{\int_0^t [E(W_n(x+X_s) - W(x+X_s))^2]\,ds\right\}$$

$$\to 0,$$

since $W_n(x+X_s) \to W(x+X_s)$, $ds \times P \times \mathbb{P}$ a.e., as $n \to \infty$, and moreover the sequence $\{(W_n(x+X_s))^2, n \geq 1\}$ is $ds \times P \times \mathbb{P}$ uniformly integrable on $[0,t] \times \Sigma \times \Omega$, thanks to (3.10).

Finally, from the occupation time formula for continuous semimartingales (see e.g. Corollary 1.6, page 209 in [15]), with $\{L_t^\cdot\}$ denoting the the local time of the process $\{X_s,\ 0 \leq s \leq t\}$,

$$\int_0^t W_n'(x+X_s)\,ds = \int_\mathbb{R} L_t^{y-x} W_n'(y)\,dy \to \int_\mathbb{R} L_t^{y-x} W(dy),$$

in $L^2(\Sigma)$, $\mathbb{P}$ a.s., as $n \to \infty$ (for more details see Section 5.7 in [10]). We used again the fact that the Brownian motions $\{X_t, t \geq 0\}$ and $\{W(y), y \in \mathbb{R}\}$ are independent. Passing now to the limit in the formula (3.11) we get the desired result. □

We can finally proceed with the following proof.

**Proof of Theorem 3.1.** Since the mapping

$$\Lambda : C(\mathbb{R}_+) \times C(\mathbb{R}) \to C(\mathbb{R}_+)$$

defined by

$$\Lambda(x, f)(t) = f(x(t))$$

is continuous, if we equip the three spaces with the topology of uniform convergence on compact sets, we first conclude from Proposition 3.7 that

$$W_\varepsilon(X_\cdot^{\varepsilon,x}) \Rightarrow W(X_\cdot^x).$$

Hence from the formulas for $\langle M^{\varepsilon,x}\rangle_t$ and $\langle M^{\varepsilon,x}, Z^\varepsilon\rangle_t$ above, and Lemma 3.9, we deduce easily that

$$\langle M^{\varepsilon,x}\rangle_\cdot \Rightarrow \frac{\overline{c}^2}{\tilde{a}} \int_0^\cdot W(X_s^x)\,ds,$$

$$\langle M^{\varepsilon,x}, Z^\varepsilon\rangle_\cdot \Rightarrow \overline{c} \int_0^\cdot W(X_s^x)\,ds,$$

and consequently

$$M_\cdot^{\varepsilon,x} \Rightarrow \frac{\overline{c}}{\tilde{a}} \int_0^\cdot W(X_s^x)\,dX_s.$$

From Proposition 3.7, those convergences are joint with those of $(X_\cdot^{\varepsilon,x}, F_\varepsilon)$. Consequently

$$F_\varepsilon(X_t^{\varepsilon,x}) - F_\varepsilon(x) - M_t^{\varepsilon,x} \Rightarrow F(X_t^x) - F(x) - \frac{\overline{c}}{\tilde{a}} \int_0^t W(X_s^x)\,dX_s.$$

The convergence $Y_t^{\varepsilon,x} \Rightarrow Y_t^x$ now follows from (3.6) and (3.8). The result finally follows from the fact that all the above convergences are joint with that of $X_\cdot^{\varepsilon,x}$. □

## 4. Pointwise convergence of the sequence $u^\varepsilon$

The first part of this section is devoted to establishing uniform integrability estimates for the exponent in the Feynman–Kac formula (Propositions 4.4 and 4.5) which are essential for the proof of the pointwise convergence part of Theorem 2.2, to which the second part of this section is devoted.

We first define the following $\mathbb{R}_+$-valued random variables, for $0 < \gamma < 1/2$, $\varepsilon \geq 0$:

$$\xi_{\gamma,\varepsilon} = \sup_{x \in \mathbb{R}} \frac{|W_\varepsilon(x)|}{(1+|x|)^{1-\gamma}}.$$

We have the following lemma.

**Lemma 4.1.** *For any $0 < \gamma < 1/2$ and $\varepsilon_0 > 0$, the collection of random variables $\{\xi_{\gamma,\varepsilon}, 0 < \varepsilon \leq \varepsilon_0\}$ is tight.*

**Proof.** Due to the symmetry it is sufficient to estimate $|W_\varepsilon(x)|$ for $x > 0$. We have

$$E(|W_\varepsilon(r)|^2) = \varepsilon \int_0^{r/\varepsilon} \int_0^{r/\varepsilon} E(c(s)c(t))\,ds\,dt$$

$$\leq 2\varepsilon \int_0^{r/\varepsilon} \int_0^\infty |E(c(0)c(s))|\,ds\,dt$$

$$\leq 2rc_0.$$

Denote

$$\eta_t = \int_0^\infty E(c(s+t)|\mathcal{G}_t)\,ds.$$

By Proposition 7.2.6. in [6] the process $\eta_t$ is stationary and $|\eta_t| \leq c_1$ a.s. with a non-random constant $c_1$. Moreover,

$$\int_0^t c(r)\,\mathrm{d}r - \eta_t$$

is a square integrable $\mathcal{G}_t$ martingale. Denote it by $\mathcal{N}_t$. Clearly,

$$W_\varepsilon(t) = \frac{\sqrt{\varepsilon}}{\bar{c}} \int_0^{t/\varepsilon} c(s)\,\mathrm{d}s = \frac{\sqrt{\varepsilon}}{\bar{c}} \mathcal{N}_{t/\varepsilon} + \frac{\sqrt{\varepsilon}}{\bar{c}} \eta_{t/\varepsilon},$$

and thus we deduce from Doob's inequality

$$E\Big(\sup_{0\leq t\leq r} |W_\varepsilon(t)|^2\Big) \leq \frac{2}{\bar{c}^2} E\Big(\sup_{0\leq t\leq r/\varepsilon} (\sqrt{\varepsilon}\mathcal{N}_t)^2\Big) + 2\frac{c_1^2 \varepsilon}{\bar{c}^2}$$

$$\leq \frac{4}{\bar{c}^2} E((\sqrt{\varepsilon}\mathcal{N}_{r/\varepsilon})^2) + 2\frac{c_1^2 \varepsilon}{\bar{c}^2}$$

$$\leq 8E(|W_\varepsilon(r)|^2) + 10\frac{c_1^2 \varepsilon}{\bar{c}^2}$$

$$\leq C(\varepsilon + r),$$

provided $C = (16c_0) \vee (10c_1^2/\bar{c}^2)$. Now for $j \geq 1$, $M > 0$,

$$P\Big(\sup_{2^{j-1}<r\leq 2^j} \frac{|W_\varepsilon(r)|}{(1+r)^{1-\gamma}} \geq M\Big) \leq P\Big(\sup_{0\leq r\leq 2^j} |W_\varepsilon(r)| \geq (1+2^{j-1})^{1-\gamma} M\Big)$$

$$\leq \frac{C(\varepsilon + 2^j)}{M^2(1+2^{j-1})^{2-2\gamma}}$$

$$\leq (\varepsilon \vee 1)\frac{2C}{M^2}(1+2^{j-1})^{2\gamma-1}.$$

Summing up over $j \geq 1$, we deduce that

$$P(\xi_{\gamma,\varepsilon} \geq M) \leq 2P\Big(\sup_{r>0} \frac{|W_\varepsilon(r)|}{(1+r)^{1-\gamma}} \geq M\Big)$$

$$\leq (\varepsilon \vee 1)\frac{4C}{M^2} \sum_{j=0}^{\infty} (1+2^j)^{2\gamma-1}$$

$$\leq (\varepsilon \vee 1)\frac{C'}{M^2}.$$

The lemma is established. $\square$

**Remark 4.2.** *We can in fact show that, as $\varepsilon \to 0$,*

$$\xi_{\gamma,\varepsilon} \Rightarrow \xi_\gamma := \sup_{x\in\mathbb{R}} \frac{|W(x)|}{(1+|x|)^{1-\gamma}},$$

*provided again $0 < \gamma < 1/2$, but we shall not use that result.*

We next state a result, which is an immediate consequence of Lemma 3.2.

**Lemma 4.3.** *There exists a continuous mapping*

$$\rho : \mathbb{R}_+ \times \mathbb{R}_+ \times (0,2) \to \mathbb{R}_+$$

*such that for all $c > 0$, $t > 0$, $\varepsilon > 0$, $0 < p < 2$,*

$$\mathbb{E}^{\varepsilon,\cdot} \exp(c|X_t^{\varepsilon,x}|^p) \le \rho(c,t,p).$$

We next establish the following proposition.

**Proposition 4.4.** *For any $0 < \gamma < 1/2$, there exists a continuous mapping $\Psi_\gamma : \mathbb{R}_+ \to \mathbb{R}_+$ such that*

$$\mathbb{E}^{\varepsilon,\cdot}\left[\left(\exp\left(\frac{1}{\sqrt{\varepsilon}}\int_0^t c\left(\frac{X_s^{\varepsilon,x}}{\varepsilon}\right)\mathrm{d}s\right)\right)^2\right] \le \Psi_\gamma(\xi_{\gamma,\varepsilon}). \quad (4.1)$$

**Proof.** Since from (3.6) and (3.5)

$$\frac{1}{\sqrt{\varepsilon}}\int_0^t c\left(\frac{X_s^{\varepsilon,x}}{\varepsilon}\right)\mathrm{d}s = 2\left[\bar{c}\int_x^{X_t^{\varepsilon,x}} \frac{W_\varepsilon(y)}{a(y/\varepsilon)}\mathrm{d}y - M_t^\varepsilon\right],$$

we obtain

$$\mathbb{E}^{\varepsilon,\cdot}\left(\exp\left(\frac{1}{\sqrt{\varepsilon}}\int_0^t c\left(\frac{X_s^{\varepsilon,x}}{\varepsilon}\right)\mathrm{d}s\right)^2\right) = \mathbb{E}^{\varepsilon,\cdot}\exp\left(4\bar{c}\int_x^{X_t^{\varepsilon,x}} \frac{W_\varepsilon(y)}{a(y/\varepsilon)}\mathrm{d}y - 4M_t^\varepsilon\right)$$

$$\le \left(\mathbb{E}^{\varepsilon,\cdot}\exp\left(8\bar{c}\int_x^{X_t^{\varepsilon,x}} \frac{W_\varepsilon(y)}{a(y/\varepsilon)}\mathrm{d}y\right)\right)^{1/2}(\mathbb{E}^{\varepsilon,\cdot}\exp(-8M_t^\varepsilon))^{1/2}. \quad (4.2)$$

Clearly, it suffices to estimate each factor on the r.h.s. of (4.2) separately.

$$\mathbb{E}^{\varepsilon,\cdot}\exp\left(8\bar{c}\int_x^{X_t^{\varepsilon,x}}\frac{W_\varepsilon(y)}{a(y/\varepsilon)}\mathrm{d}y\right) \le \mathbb{E}^{\varepsilon,\cdot}\exp\left(8\frac{\bar{c}}{\underline{c}}\int_x^{X_t^{\varepsilon,x}}|W_\varepsilon(y)|\,\mathrm{d}y\right)$$

$$\le \mathbb{E}^{\varepsilon,\cdot}\exp\left(8\frac{\bar{c}}{\underline{c}}\int_x^{X_t^{\varepsilon,x}}\xi_{\gamma,\varepsilon}(1+|y|)^{1-\gamma}\,\mathrm{d}y\right)$$

$$\le \mathbb{E}^{\varepsilon,\cdot}\exp\left(c'\frac{\xi_{\gamma,\varepsilon}}{2-\gamma}[(1+|X_t^{\varepsilon,x}|)^{2-\gamma} - (1+|x|)^{2-\gamma}]\right)$$

$$\le \Psi_\gamma^1(\xi_{\gamma,\varepsilon}),$$

where we have used Lemma 3.2 for the last inequality. The second factor on the r.h.s. of (4.2) can be estimated as follows

$$\mathbb{E}^{\varepsilon,\cdot}\exp(-8M_t^\varepsilon) \le (\mathbb{E}^{\varepsilon,\cdot}\exp(-16M_t^\varepsilon - 128\langle M^\varepsilon\rangle_t))^{1/2}(\mathbb{E}^{\varepsilon,\cdot}\exp(128\langle M^\varepsilon\rangle_t))^{1/2}.$$

The first term on the r.h.s. does not exceed 1. For the second one we have by Jensen's inequality

$$\mathbb{E}^{\varepsilon,\cdot}\exp(128\langle M^\varepsilon\rangle_t) \le \mathbb{E}^{\varepsilon,\cdot}\exp\left(c''\int_0^t W_\varepsilon^2(X_s^{\varepsilon,x})\,\mathrm{d}s\right)$$

$$\le t^{-1}\int_0^t \mathbb{E}^{\varepsilon,\cdot}\exp(c''t W_\varepsilon^2(X_s^{\varepsilon,x}))\,\mathrm{d}s$$

$$\le \sup_{0\le s\le t}\mathbb{E}^{\varepsilon,\cdot}\exp(c''t[\xi_{\gamma,\varepsilon}]^2(1+|X_s^{\varepsilon,x}|)^{2-2\gamma})$$

$$\le \Psi_\gamma^2(\xi_{\gamma,\varepsilon}),$$

where we have again used Lemma 3.2 for the last inequality. The result clearly follows. □

Clearly, the same proof allows us to establish the slightly more general proposition.

**Proposition 4.5.** *Let $\{\tau^\varepsilon, \varepsilon \geq 0\}$ be a collection of stopping times such that $0 \leq \tau^\varepsilon \leq t$ $P \times \mathbb{P}$ a.s. Then for any $0 < \gamma < 1/2$,*

$$\mathbb{E}^{\varepsilon,\cdot}\left[\left(\exp\left(\frac{1}{\sqrt{\varepsilon}}\int_0^{\tau^\varepsilon} c\left(\frac{X_s^{\varepsilon,x}}{\varepsilon}\right)\mathrm{d}s\right)\right)^2\right] \leq \Psi_\gamma(\xi_{\gamma,\varepsilon}),$$

*where $\Psi_\gamma : \mathbb{R}_+ \to \mathbb{R}_+$ is the mapping which appeared in (4.1).*

We can now proceed with the following proof.

**Proof of the pointwise convergence in Theorem 2.2.** We will now show that for each $(t,x) \in \mathbb{R}_+ \times \mathbb{R}$, $u^\varepsilon(t,x) \Rightarrow u(t,x)$, as $\varepsilon \to 0$. We delete the parameters $t$ and $x$. It suffices to show that for any $\varphi \in C(\mathbb{R};[0,1])$, $\varphi$ Lipschitz continuous, as $\varepsilon \to 0$,

$$E\varphi(\mathbb{E}^{\varepsilon,\cdot}[g(X^\varepsilon)\exp(Y^\varepsilon)]) \to E\varphi(\mathbb{E}^{\varepsilon,\cdot}[g(X)\exp(Y)]), \tag{4.3}$$

where $X^\varepsilon = X_t^{\varepsilon,x}$, $X = X_t^x = x + \sqrt{\tilde{a}}B_t$, and (recall (3.6) and (3.5))

$$Y^\varepsilon = 2\left[\overline{c}\int_x^{X^\varepsilon}\frac{W_\varepsilon(y)}{a(y/\varepsilon)}\mathrm{d}y - M_t^\varepsilon\right] = 2\overline{c}\left(\int_x^{X_t^{\varepsilon,x}}\frac{W_\varepsilon(y)}{a(y/\varepsilon)}\mathrm{d}y - \frac{1}{\tilde{a}}\int_0^t W_\varepsilon(X_s^{\varepsilon,x})\mathrm{d}Z_s^\varepsilon\right),$$

$$Y = 2\frac{\overline{c}}{\tilde{a}}\left[\int_x^X W(y)\mathrm{d}y - \int_0^t W(x+X_s)\mathrm{d}X_s\right].$$

The fact that this $Y$ equals the exponent in the Feynman–Kac formula for $u(t,x)$ follows from (3.8). We first approximate $Y^\varepsilon$ by $Y^{\varepsilon,M}$ as follows. For each $\varepsilon > 0$, $M > 0$, let

$$\tau_M^\varepsilon = \inf\{s \geq 0; |X_s^{\varepsilon,x}| \geq M\}$$

and

$$Y^{\varepsilon,M} = 2\overline{c}\left(\int_x^{X_{t\wedge\tau_M^\varepsilon}^{\varepsilon,x}}\frac{W_\varepsilon(y)}{a(y/\varepsilon)}\mathrm{d}y - \frac{1}{\tilde{a}}\int_0^{t\wedge\tau_M^\varepsilon} W_\varepsilon(X_s^{\varepsilon,x})\mathrm{d}Z_s^\varepsilon\right).$$

We postpone the proof of the following lemma.

**Lemma 4.6.**

$$\sup_{\varepsilon>0}|E\varphi(\mathbb{E}^{\varepsilon,\cdot}[g(X_t^{\varepsilon,x})\exp(Y^\varepsilon)]) - E\varphi(\mathbb{E}^{\varepsilon,\cdot}[g(X_t^{\varepsilon,x})\exp(Y^{\varepsilon,M})])| \to 0,$$

*as $M \to \infty$.*

Since the collection of random processes $\{W_\varepsilon(y); y \in \mathbb{R}\}$ is $P$-tight, for all $\delta > 0$, there exists $N \in \mathbb{N}$ and $f_{\delta,1}, f_{\delta,2}, \ldots, f_{\delta,N} \in C_b(\mathbb{R}_+)$ such that if

$$\tilde{B}_k^{\delta,\varepsilon} := \Big\{\sup_{-\delta^{-1}\leq x\leq\delta^{-1}}|W_\varepsilon(x) - f_{\delta,k}(x)| \leq \delta\Big\}, \quad 1 \leq k \leq N,$$

then

$$P\left(\left(\bigcup_{k=1}^N \tilde{B}_k^{\delta,\varepsilon}\right)^c\right) \leq \delta.$$

Let now
$$B_1^{\delta,\varepsilon} = \tilde{B}_1^{\delta,\varepsilon}$$

and for $2 \le k \le N$,
$$B_k^{\delta,\varepsilon} = \tilde{B}_k^{\delta,\varepsilon} \setminus \bigcup_{i=1}^{k-1} B_i^{\delta,\varepsilon},$$

and finally
$$A^{\delta,\varepsilon} := \left( \bigcup_{k=1}^{N} B_k^{\delta,\varepsilon} \right)^c,$$

so that $P(A^{\delta,\varepsilon}) \le \delta$. Note that $\|f_{\delta,k}\|_\infty$ depends on $\delta$. However, we can and will assume that for some $0 < \gamma < 1/2$,

$$|f_{\delta,k}(x)| \le 2\xi_\gamma (1+|x|)^{1-\gamma}, \tag{4.4}$$

for all $\delta > 0$, $k \in \mathbb{N}$.

We now develop
$$E\varphi(\mathbb{E}^{\varepsilon,\cdot}[g(X^\varepsilon)e^{Y^{\varepsilon,\delta^{-1}}}]) = \sum_{k=1}^{N} E\{\varphi(\mathbb{E}^{\varepsilon,\cdot}[g(X^\varepsilon)e^{Y^{\varepsilon,\delta^{-1}}}]); B_k^{\delta,\varepsilon}\} + E\{\varphi(\mathbb{E}^{\varepsilon,\cdot}[g(X^\varepsilon)e^{Y^{\varepsilon,\delta^{-1}}}]); A^{\delta,\varepsilon}\}.$$

The last term in the above right-hand side is bounded in absolute value by $\delta$. Now for $1 \le k \le N$,

$$E\{\varphi(\mathbb{E}^{\varepsilon,\cdot}[g(X^\varepsilon)e^{Y^{\varepsilon,\delta^{-1}}}]); B_k^{\delta,\varepsilon}\} = E\{\varphi(\mathbb{E}^{\varepsilon,\cdot}[g(X^\varepsilon)e^{Y_k^{\delta,\varepsilon}}]); B_k^{\delta,\varepsilon}\} + e_k^{\varepsilon,\delta},$$

where
$$Y_k^{\delta,\varepsilon} := 2\overline{c}\left( \int_x^{X_t^{\varepsilon,x}} \frac{f_{\delta,k}(y)}{a(y/\varepsilon)}\,dy - \frac{1}{\tilde{a}} \int_0^t f_{\delta,k}(X_s^{\varepsilon,x})\,dZ_s^\varepsilon \right).$$

We postpone the proofs of following lemmas.

**Lemma 4.7.** *There exists a constant $C$, which depends only on $t$, $\|f_{\delta,k}\|_\infty$ and the constants appearing in (2.1), such that for each $\delta > 0$, $1 \le k \le N$, $\sigma \in \Sigma$,*

$$\sup_{\varepsilon > 0} \mathbb{E}^{\varepsilon,\sigma}[e^{2Y_k^{\delta,\varepsilon}}] \le C, \quad P \text{ a.s.}$$

**Lemma 4.8.**
$$\sup_{\varepsilon \ge 0} \left| \sum_{k=1}^{N} e_k^{\varepsilon,\delta} \right| \to 0,$$

*as $\delta \to 0$.*

It follows readily from Propositions 3.10 and 3.11 together with Lemma 4.7 that, as $\varepsilon \to 0$,

$$\mathbb{E}^{\varepsilon,\sigma}[g(X^\varepsilon)e^{Y_k^{\delta,\varepsilon}}] \to \mathbb{E}^{\varepsilon,\sigma}[g(X)e^{Y_k^\delta}],$$

$P$ a.s., where

$$Y_k^\delta = 2\frac{\overline{c}}{\widetilde{a}}\left[\int_x^{X_t^x} f_{\delta,k}(y)\,\mathrm{d}y - \int_0^t f_{\delta,k}(X_s^x)\,\mathrm{d}X_s^x\right].$$

Let $B_k^\delta$, $1 \leq k \leq N$, denote the sets defined exactly as the $B_k^{\delta,\varepsilon}$'s, but with $W_\varepsilon$ replaced by $W$. The boundaries of those sets being of zero Wiener measure, we conclude from the last statement and the fact that $W_\varepsilon \Rightarrow W$ that as $\varepsilon \to 0$,

$$E\{\varphi(\mathbb{E}^{\varepsilon,\cdot}[g(X^\varepsilon)\mathrm{e}^{Y_k^{\delta,\varepsilon}}]); B_k^{\delta,\varepsilon}\} \to E\{\varphi(\mathbb{E}^{\varepsilon,\cdot}[g(X)\mathrm{e}^{Y_k^\delta}]); B_k^\delta\}.$$

Now, in the same way as above we obtain

$$E\varphi(\mathbb{E}^{\varepsilon,\cdot}[g(X)\mathrm{e}^Y]) = \sum_{k=1}^N E\{\varphi(\mathbb{E}^{\varepsilon,\cdot}[g(X)\mathrm{e}^Y]); B_k^\delta\} + E\{\varphi(\mathbb{E}^{\varepsilon,\cdot}[g(X)\mathrm{e}^Y]); A^\delta\},$$

$P(A^\delta) \leq \delta$, and for each $1 \leq k \leq N$,

$$E\{\varphi(\mathbb{E}^{\varepsilon,\cdot}[g(X)\mathrm{e}^Y]); B_k^\delta\} = E\{\varphi(\mathbb{E}^{\varepsilon,\cdot}[g(X)\mathrm{e}^{Y_k^\delta}]); B_k^\delta\} + e_k^\delta.$$

All we need to conclude the proof is the next lemma.

**Lemma 4.9.**

$$\left|\sum_{k=1}^N e_k^\delta\right| \to 0,$$

as $\delta \to 0$.

It remains to prove the four lemmas.

**Proof of Lemma 4.6.** Let $K$ denote the product of the Lipschitz constant of $\varphi$ by $\|g\|_\infty$. Since $\varphi$ takes values in $[0,1]$,

$$E|\varphi(\mathbb{E}^{\varepsilon,\cdot}[g(X^\varepsilon)\mathrm{e}^{Y^\varepsilon}]) - \varphi(\mathbb{E}^{\varepsilon,\cdot}[g(X^\varepsilon)\mathrm{e}^{Y^{\varepsilon,M}}])|$$

$$\leq E\{1 \wedge K\mathbb{E}^{\varepsilon,\cdot}(|\mathrm{e}^{Y^\varepsilon} - \mathrm{e}^{Y^{\varepsilon,M}}|)\}$$

$$\leq E\{1 \wedge \sqrt{2}K(\mathbb{E}^{\varepsilon,\cdot}[\mathrm{e}^{2Y^\varepsilon} + \mathrm{e}^{2Y^{\varepsilon,M}}])^{1/2}(\mathbb{P}^{\varepsilon,\cdot}(\tau_M^\varepsilon < t))^{1/2}\}$$

$$\leq E\{1 \wedge [2K(\Psi_\gamma(\xi_{\gamma,\varepsilon}))^{1/2}(\mathbb{P}^{\varepsilon,\cdot}(\tau_M^\varepsilon < t))^{1/2}]\}$$

$$\leq L(\mathbb{P}^{\varepsilon,\cdot}(\tau_M^\varepsilon < t))^{1/2} + P(\Psi_\gamma(\xi_{\gamma,\varepsilon}) > L^2/4K^2),$$

where $0 < \gamma < 1/2$ and $\Psi_\gamma(\xi_{\gamma,\varepsilon})$ appears in Propositions 4.4 and 4.5. Since the latter is $P$-tight (for $\varepsilon > 0$), we can choose $L = L_\eta$ such that the second term in the last expression above is less than $\eta/2$. It remains to choose $M$ large enough (exploiting this time the tightness of $\{X^{\varepsilon,x}; \varepsilon > 0\}$) such that

$$L_\eta(\mathbb{P}^{\varepsilon,\cdot}(\tau_M^\varepsilon < t))^{1/2} \leq \eta/2. \qquad \square$$

**Proof of Lemma 4.7.** There exist two constants $c_1$ and $c_2$ such that

$$\mathbb{E}^{\varepsilon,\cdot}\mathrm{e}^{2Y_k^{\delta,\varepsilon}} \leq (\mathbb{E}^{\varepsilon,\cdot}\mathrm{e}^{c_1|X_t^{\varepsilon,x}-x|})^{1/2}\left(\mathbb{E}^{\varepsilon,\cdot}\mathrm{e}^{c_2\int_0^t f_{\delta,k}(X_s^{\varepsilon,x})\,\mathrm{d}Z_s^\varepsilon}\right)^{1/2}.$$

The bound for the first factor on the right follows from Lemma 3.2, and the bound for the second factor follows easily from the boundedness of both $f_{\delta,k}$ and $\frac{d}{ds}\langle Z^\varepsilon\rangle_s$. $\square$

**Proof of Lemma 4.8.** We have, by an argument similar to that in the proof of Lemma 4.6,

$$\sum_{k=1}^N |e_k^{\varepsilon,\delta}| \leq E\left(1 \wedge \sum_{k=1}^N K\mathbb{E}^{\varepsilon,\cdot}[|e^{Y^{\varepsilon,\delta^{-1}}} - e^{Y_k^{\delta,\varepsilon}}|]\mathbf{1}_{B_k^{\delta,\varepsilon}}\right)$$

$$\leq E\left(1 \wedge \sum_{k=1}^N K\mathbb{E}^{\varepsilon,\cdot}[(e^{Y^{\varepsilon,\delta^{-1}}} + e^{Y_k^{\delta,\varepsilon}})|Y^{\varepsilon,M} - Y_k^{\delta,\varepsilon}|]\mathbf{1}_{B_k^{\delta,\varepsilon}}\right).$$

Now

$$\mathbb{E}^{\varepsilon,\cdot}[(e^{Y^{\varepsilon,\delta^{-1}}} + e^{Y_k^{\delta,\varepsilon}})|Y^{\varepsilon,\delta^{-1}} - Y_k^{\delta,\varepsilon}|]\mathbf{1}_{B_k^{\delta,\varepsilon}}$$

$$\leq (2\mathbb{E}^{\varepsilon,\cdot}[e^{2Y^{\varepsilon,\delta^{-1}}} + e^{2Y_k^{\delta,\varepsilon}}])^{1/2}(\mathbf{1}_{B_k^{\delta,\varepsilon}}\mathbb{E}^{\varepsilon,\cdot}[|Y^{\varepsilon,\delta^{-1}} - Y_k^{\delta,\varepsilon}|^2])^{1/2}$$

$$\leq c(2[\Psi_\gamma(\xi_{\gamma,\varepsilon}) + C])^{1/2}$$

$$\times (\delta^2\mathbb{E}^{\varepsilon,\cdot}(|X_{t\wedge\tau_{\delta^{-1}}^\varepsilon}^{\varepsilon,x}|^2) + \mathbb{E}^{\varepsilon,\cdot}(|X_t^{\varepsilon,x} - X_{t\wedge\tau_{\delta^{-1}}^\varepsilon}^{\varepsilon,x}|^2) + \delta t + \mathbb{E}^{\varepsilon,\cdot}(t - t\wedge\tau_{\delta^{-1}}^\varepsilon))^{1/2}\mathbf{1}_{B_k^{\delta,\varepsilon}}.$$

Finally,

$$\sum_{k=1}^N |e_k^{\varepsilon,\delta}| \leq E(1 \wedge [K(1 + \Psi_\gamma(\xi_{\gamma,\varepsilon}))^{1/2}(\delta t + C\delta^2 + \rho(t,\delta))^{1/2}]),$$

where

$$\rho(t,\delta) = \sup_{\varepsilon>0}[\mathbb{E}^{\varepsilon,\cdot}(|X_t^{\varepsilon,x} - X_{t\wedge\tau_{\delta^{-1}}^\varepsilon}^{\varepsilon,x}|^2) + \mathbb{E}^{\varepsilon,\cdot}(t - t\wedge\tau_{\delta^{-1}}^\varepsilon)]$$

$$\to 0,$$

as $\delta \to 0$. The end of the proof is similar to that of Lemma 4.6. $\square$

**Proof of Lemma 4.9.** This proof is similar to that of Lemma 4.8.

$$\sum_{k=1}^N E\{|\varphi(\mathbb{E}^{\varepsilon,\cdot}[g(X)e^Y]) - \varphi(\mathbb{E}^{\varepsilon,\cdot}[g(X)e^{Y_k^\delta}])|; B_k^\delta\} \leq E\left\{1 \wedge K\sum_{k=1}^N \mathbb{E}^{\varepsilon,\cdot}|e^Y - e^{Y_k^\delta}|\mathbf{1}_{B_k^\delta}\right\}.$$

Now

$$\mathbf{1}_{B_k^\delta}\mathbb{E}^{\varepsilon,\cdot}|e^Y - e^{Y_k^\delta}| \leq (2\mathbb{E}^{\varepsilon,\cdot}(e^{2Y} + e^{2Y_k^\delta}))^{1/2}(\mathbf{1}_{B_k^\delta}\mathbb{E}^{\varepsilon,\cdot}[|Y - Y_k^\delta|^2])^{1/2},$$

and on the set $B_k^\delta$, using in particular (4.4),

$$\mathbb{E}^{\varepsilon,\cdot}(|Y - Y_k^\delta|^2) \leq 8\left(\frac{\overline{c}}{\widetilde{a}}\right)^2 \mathbb{E}^{\varepsilon,\cdot}\left[\left(\int_x^{x+B_t}[W(y) - f_{\delta,k}(y)]\,dy\right)^2\right]$$

$$+ 8\left(\frac{\overline{c}}{\widetilde{a}}\right)^2 \mathbb{E}^{\varepsilon,\cdot}\int_0^t |W(x + B_s) - f_{\delta,k}(x + B_s)|^2\,ds$$

$$\leq c\Big(\delta^2 t + \xi_\gamma t^{1/2} \vee t^{1-\gamma/2}\sqrt{\mathbb{P}^{\varepsilon,\cdot}(|B_t| > \delta^{-1} - |x|)}$$

$$+ c\xi_\gamma^2 t^{2-\gamma}\mathbb{P}^{\varepsilon,\cdot}\left(\sup_{0\leq s\leq t}|x + B_s| > \delta^{-1}\right)\Big),$$

which clearly goes to 0, as $\delta \to 0$. The result follows. $\square$

$\square$

## 5. Convergence in $C(\mathbb{R}_+ \times \mathbb{R})$

It remains both to prove the convergence of the finite dimensional distributions of $u^\varepsilon$ towards those of $u$, and to establish that the sequence $\{u^\varepsilon; \varepsilon > 0\}$ is tight as a collection of random elements of $C(\mathbb{R}_+ \times \mathbb{R})$.

**Theorem 5.1.** *For any $\ell \in \mathbb{N}$, $(t_i, x_i) \in \mathbb{R}_+ \times \mathbb{R}$, $1 \leq i \leq \ell$,*

$$(u^\varepsilon(t_1, x_1), \ldots, u^\varepsilon(t_\ell, x_\ell)) \Rightarrow (u(t_1, x_1), \ldots, u(t_\ell, x_\ell))$$

*as $\varepsilon \to 0$.*

**Proof.** We only sketch the proof, the details being identical to those of the proof of the pointwise convergence, as were given in the previous section. For each $1 \leq i \leq \ell$, we define $X_i^\varepsilon := X_{t_i}^{\varepsilon, x_i}$ and

$$Y_i^\varepsilon := \frac{1}{\sqrt{\varepsilon}} \int_0^{t_i} c\left(\frac{X_s^{\varepsilon, x_i}}{\varepsilon}\right) ds.$$

We need to take the limit as $\varepsilon \to 0$ in the quantity

$$E\varphi(\mathbb{E}^{\varepsilon, \cdot}[g(X_1^\varepsilon) \exp(Y_1^\varepsilon)], \ldots, \mathbb{E}^{\varepsilon, \cdot}[g(X_\ell^\varepsilon) \exp(Y_\ell^\varepsilon)]),$$

where $\varphi \in C(\mathbb{R}^\ell; [0, 1])$ is Lipschitz continuous. For that sake, referring to the notations in the previous section, for each $\delta > 0$, $1 \leq k \leq N$, we define

$$Y_{i,k}^{\delta, \varepsilon} := 2\bar{c}\left(\int_{x_i}^{X_i^\varepsilon} \frac{f_{\delta, k}(z)}{a(z/\varepsilon)} dz - \frac{1}{\tilde{a}} \int_0^{t_i} f_{\delta, k}(X_s^{\varepsilon, x_i}) dZ_s^{\varepsilon, i}\right),$$

where

$$Z_t^{\varepsilon, i} = X_t^{\varepsilon, x_i} + \varepsilon \chi\left(\frac{X_t^{\varepsilon, x_i}}{\varepsilon}\right).$$

We have that for each $1 \leq i \leq \ell$,

$$\mathbb{E}^{\varepsilon, \cdot}[g(X_i^\varepsilon) \exp(Y_i^\varepsilon)] \simeq \sum_{k=1}^N \mathbb{E}^{\varepsilon, \cdot}[g(X_i^\varepsilon) \exp(Y_{i,k}^{\delta, \varepsilon})] \mathbf{1}_{B_k^\delta},$$

and consequently

$$E\varphi(\mathbb{E}^{\varepsilon, \cdot}[g(X_1^\varepsilon) \exp(Y_1^\varepsilon)], \ldots, \mathbb{E}^{\varepsilon, \cdot}[g(X_\ell^\varepsilon) \exp(Y_\ell^\varepsilon)])$$

$$= \sum_{k=1}^N E[\varphi(\mathbb{E}^{\varepsilon, \cdot}[g(X_1^\varepsilon) \exp(Y_{1,k}^{\delta, \varepsilon})], \ldots, \mathbb{E}^{\varepsilon, \cdot}[g(X_\ell^\varepsilon) \exp(Y_{\ell,k}^{\delta, \varepsilon})]) \mathbf{1}_{B_k^\delta}].$$

But for each $1 \leq i \leq \ell$, $1 \leq k \leq N$,

$$\mathbb{E}^{\varepsilon, \cdot}[g(X_i^\varepsilon) \exp(Y_{i,k}^{\delta, \varepsilon})] \to \mathbb{E}^{\varepsilon, \cdot}[g(X_i) \exp(Y_{i,k}^{\delta})]$$

$P$ a.s., as $\varepsilon \to 0$, where $X_i = x_i + \sqrt{\tilde{a}} B_{t_i}$ and

$$Y_{i,k}^{\delta} = 2\frac{\bar{c}}{\tilde{a}}\left(\int_{x_i}^{X_i} f_{\delta, k}(z) dz - \int_0^{t_i} f_{\delta, k}(x_i + \sqrt{\tilde{a}} B_s) dB_s\right).$$

The result follows. □

We finally establish the tightness result. All we need to show is the following theorem.

**Theorem 5.2.** *There exists a sequence $\rho_\varepsilon$ in $\mathbb{R}_+$ which tends to 0 as $\varepsilon \to 0$, a mapping*

$$\varphi: \mathbb{R}_+^4 \times \mathbb{R}^2 \to \mathbb{R}_+$$

*such that for all $T > 0$, $M > 0$,*

$$\sup_{0 \le s,t \le T; |x|,|y| \le M; |t-s| \le \delta; |x-y| \le \delta} \varphi(T, M, s, t, x, y) \to 0 \tag{5.1}$$

*as $\delta \to 0$, and a continuous mapping $\Phi: \mathbb{R}_+^3 \to \mathbb{R}_+$ such that for all $0 \le s, t \le T$, $|x|, |y| \le M$,*

$$|u(t,x) - u(s,y)| \le \Phi(T, M, \xi_{\gamma,\varepsilon})[\rho_\varepsilon + \varphi(T, M, s, t, x, y)]. \tag{5.2}$$

**Proof.** We will establish (5.2) with

$$\varphi(T, M, s, t, x, y) = \sqrt{|t-s|} + \sqrt{|x-y|} + \sqrt{\mathbb{E}([g(x + \sqrt{\tilde{a}}B_t) - g(y + \sqrt{\tilde{a}}B_s)]^2)},$$

which satisfies (5.1).

We need to consider

$$\begin{aligned}
u^\varepsilon(t,x) - u^\varepsilon(s,y) &= \mathbb{E}^{\varepsilon,\cdot}\left[g(X_t^{\varepsilon,x}) e^{1/\sqrt{\varepsilon} \int_0^t c(X_r^{\varepsilon,x}/\varepsilon)\,dr}\right] - \mathbb{E}^{\varepsilon,\cdot}\left[g(X_s^{\varepsilon,y}) e^{1/\sqrt{\varepsilon} \int_0^s c(X_r^{\varepsilon,y}/\varepsilon)\,dr}\right] \\
&= \mathbb{E}^{\varepsilon,\cdot}\left\{[g(X_t^{\varepsilon,x}) - g(X_s^{\varepsilon,y})] e^{1/\sqrt{\varepsilon} \int_0^t c(X_r^{\varepsilon,x}/\varepsilon)\,dr}\right\} \\
&\quad + \mathbb{E}^{\varepsilon,\cdot}\left\{g(X_s^{\varepsilon,y}) \left[e^{1/\sqrt{\varepsilon} \int_0^t c(X_r^{\varepsilon,x}/\varepsilon)\,dr} - e^{1/\sqrt{\varepsilon} \int_0^s c(X_r^{\varepsilon,y}/\varepsilon)\,dr}\right]\right\}.
\end{aligned} \tag{5.3}$$

The absolute value of the first term of the right-hand side of (5.3) is dominated by

$$\sqrt{\mathbb{E}^{\varepsilon,\cdot} e^{2/\sqrt{\varepsilon} \int_0^t c(X_r^{\varepsilon,x}/\varepsilon)\,dr}} \times \sqrt{\mathbb{E}^{\varepsilon,\cdot}([g(X_t^{\varepsilon,x}) - g(X_s^{\varepsilon,y})]^2)}.$$

The first factor in the last expression contributes to the coefficient $\Phi(T, M, \xi_{\gamma,\varepsilon})$ in (5.2), while the difference between the second factor and

$$\sqrt{\mathbb{E}^{\varepsilon,\cdot}([g(x + \sqrt{\tilde{a}}B_t) - g(y + \sqrt{\tilde{a}}B_s)]^2)}$$

contributes to $\rho_\varepsilon$. Now the absolute value of the second term in the right-hand side of (5.3) is dominated by $c\mathbb{E}(A_\varepsilon + B_\varepsilon)$, where

$$A_\varepsilon := \left|e^{1/\sqrt{\varepsilon} \int_0^t c(X_r^{\varepsilon,x}/\varepsilon)\,dr} - e^{1/\sqrt{\varepsilon} \int_0^t c(X_r^{\varepsilon,y}/\varepsilon)\,dr}\right|$$

and

$$B_\varepsilon := \left|e^{1/\sqrt{\varepsilon} \int_0^t c(X_r^{\varepsilon,y}/\varepsilon)\,dr} - e^{1/\sqrt{\varepsilon} \int_0^s c(X_r^{\varepsilon,y}/\varepsilon)\,dr}\right|.$$

Below we will use repeatedly the elementary inequality $|e^a - e^b| \le |a - b|e^{a \vee b}$. In particular, we deduce that (assuming w.l.o.g. that $0 \le s \le t$)

$$B_\varepsilon \le \left|\frac{1}{\sqrt{\varepsilon}} \int_s^t c\left(\frac{X_r^{\varepsilon,y}}{\varepsilon}\right) dr\right| \left(e^{1/\sqrt{\varepsilon} \int_0^t c(X_r^{\varepsilon,y}/\varepsilon)\,dr} \vee e^{1/\sqrt{\varepsilon} \int_0^s c(X_r^{\varepsilon,y}/\varepsilon)\,dr}\right).$$

Consequently $\mathbb{E}B_\varepsilon$ is dominated by a factor which contributes to $\Phi(T, M, \xi_{\gamma,\varepsilon})$ times the square root of

$$\mathbb{E}^{\varepsilon,\cdot}\left(\left[2\overline{c}\int_{X_s^{\varepsilon,y}}^{X_t^{\varepsilon,y}}\frac{W_\varepsilon(z)}{a(z/\varepsilon)}\,\mathrm{d}z + M_t^{\varepsilon,y} - M_s^{\varepsilon,y}\right]^2\right)$$

$$\leq c\mathbb{E}^{\varepsilon,\cdot}\left(\sup_{|z|\leq |X_t^{\varepsilon,y}|\vee|X_s^{\varepsilon,y}|}W^2(z)|X_t^{\varepsilon,y} - X_s^{\varepsilon,y}|^2 + \int_s^t W_\varepsilon^2(X_r^{\varepsilon,y})\,\mathrm{d}r\right)$$

$$\leq c\xi_{\gamma,\varepsilon}^2\mathbb{E}^{\varepsilon,\cdot}\left(\left(1 + \sup_{s\leq r\leq t}|X_r^{\varepsilon,y}|\right)^{2-2\gamma}|X_t^{\varepsilon,y} - X_s^{\varepsilon,y}|^2 + \int_s^t (1 + |X_r^{\varepsilon,y}|)^{2-2\gamma}\,\mathrm{d}r\right)$$

$$\leq c\xi_{\gamma,\varepsilon}^2|t-s|.$$

By the same argument as above, $\mathbb{E}^{\varepsilon,\cdot}A_\varepsilon$ is dominated by a factor which contributes to $\Phi(T, M, \xi_{\gamma,\varepsilon})$ times the square root of

$$\mathbb{E}^{\varepsilon,\cdot}\left(\left|\int_x^{X_t^{\varepsilon,x}}\frac{W_\varepsilon(z)}{a(z/\varepsilon)}\,\mathrm{d}z - \int_y^{X_t^{\varepsilon,y}}\frac{W_\varepsilon(z)}{a(z/\varepsilon)}\,\mathrm{d}z\right|^2 + |M_t^{\varepsilon,x} - M_t^{\varepsilon,y}|^2\right). \tag{5.4}$$

While the two laws of $\{X_\cdot^{\varepsilon,x}\}$ and $\{X_\cdot^{\varepsilon,y}\}$ are given to us, the coupling between these two processes is at our disposal. We make the following choice. Given $\sigma \in \Sigma$, the two processes $\{X_\cdot^{\varepsilon,x}\}$ and $\{X_\cdot^{\varepsilon,y}\}$ are mutually independent, until the first time $\tau_{xy}^\varepsilon$ when they meet, and then the two processes follow the same trajectory. We note that as $\varepsilon \to 0$,

$$\tau_{xy}^\varepsilon \Rightarrow \tau_{xy} = \inf\{r; x + \sqrt{\tilde{a}}B_r^1 = y + \sqrt{\tilde{a}}B_r^2\},$$

where $\{B^1\}$ and $\{B^2\}$ are two mutually independent standard Brownian motions. This follows from the fact that the first time when a two-dimensional Brownian motion meets the diagonal of $\mathbb{R}^2$ is a.s. a continuous function of the trajectory. Suppose w.l.o.g. that $y \leq x$. Then

$$\mathbb{P}^{\varepsilon,\cdot}(\tau_{xy} > \delta) = \mathbb{P}^{\varepsilon,\cdot}\left(\sup_{0\leq r\leq \delta}(B_r^2 - B_r^1) < \frac{x-y}{\sqrt{\tilde{a}}}\right)$$

$$= 1 - \mathbb{P}^{\varepsilon,\cdot}\left(\sup_{0\leq r\leq \delta}(B_r^2 - B_r^1) \geq \frac{x-y}{\sqrt{\tilde{a}}}\right)$$

$$= 1 - 2\mathbb{P}^{\varepsilon,\cdot}\left(B_\delta^2 - B_\delta^1 \geq \frac{x-y}{\sqrt{\tilde{a}}}\right)$$

$$= \mathbb{P}^{\varepsilon,\cdot}\left(|Z| \leq \frac{x-y}{\sqrt{2\delta\tilde{a}}}\right)$$

$$\leq \frac{1}{\sqrt{\pi\tilde{a}}} \times \frac{|x-y|}{\sqrt{\delta}},$$

where $Z$ is an $N(0,1)$ r.v. On the other hand,

$$\mathbb{E}^{\varepsilon,\cdot}(\tau_{xy} \wedge t) \leq c\sqrt{t}|x-y|.$$

We now estimate the first term in (5.4) for $|x|, |y| \leq M$, $0 \leq t \leq T$,

$$\mathbb{E}^{\varepsilon,\cdot}\left(\left|\int_x^{X_t^{\varepsilon,x}}\frac{W_\varepsilon(z)}{a(z/\varepsilon)}\,\mathrm{d}z - \int_y^{X_t^{\varepsilon,y}}\frac{W_\varepsilon(z)}{a(z/\varepsilon)}\,\mathrm{d}z\right|^2\right)$$

$$\leq \left|\int_x^y \frac{W_\varepsilon(z)}{a(z/\varepsilon)}\,\mathrm{d}z\right|^2 + \mathbb{E}^{\varepsilon,\cdot}\left(\left|\int_x^{X_t^{\varepsilon,x}}\frac{W_\varepsilon(z)}{a(z/\varepsilon)}\,\mathrm{d}z + \int_y^{X_t^{\varepsilon,y}}\frac{W_\varepsilon(z)}{a(z/\varepsilon)}\,\mathrm{d}z\right|^2 \mathbf{1}_{\{\tau_{xy}^\varepsilon > t\}}\right)$$

$$\leq c(M)\xi_{\gamma,\varepsilon}^2(|x-y|^2 + \mathbb{E}^{\varepsilon,\cdot}[(1+|x|\vee|X_t^{\varepsilon,x}|\vee|y|\vee|X_t^{\varepsilon,y}|)^{4-2\gamma}\mathbf{1}_{\{\tau_{xy}^\varepsilon>t\}}])$$
$$\leq c(M,T)\xi_{\gamma,\varepsilon}^2(|x-y|^2 + |x-y| + \rho'_\varepsilon),$$

where $\rho'_\varepsilon \to 0$, as $\varepsilon \to 0$. We finally estimate the second term in (5.4), again for $|x|, |y| \leq M$, $0 \leq t \leq T$,

$$\mathbb{E}^{\varepsilon,\cdot}(|M_t^{\varepsilon,x} - M_t^{\varepsilon,y}|^2) \leq c\mathbb{E}^{\varepsilon,\cdot}\int_0^{\tau_{xy}^\varepsilon \wedge t}[W_\varepsilon^2(X_s^{\varepsilon,x}) + W_\varepsilon^2(X_s^{\varepsilon,y})]\,\mathrm{d}s$$
$$\leq c\xi_{\gamma,\varepsilon}^2\mathbb{E}^{\varepsilon,\cdot}\left[\sup_{0\leq s\leq t}(1+|X_s^{\varepsilon,x}|+|X_s^{\varepsilon,y}|)^{2-2\gamma}(\tau_{xy}^\varepsilon \wedge t)\right]$$
$$\leq c(M,T)\xi_{\gamma,\varepsilon}^2(|x-y| + \rho'_\varepsilon).$$

The theorem is established. □

## 6. The stochastic PDE for the limit $u$

In this section, we study the limiting SPDE. Formally, we would expect that it reads

$$\frac{\partial u}{\partial t}(t,x) = \frac{\tilde{a}}{2}\frac{\partial^2 u}{\partial x^2}(t,x) + \overline{c}u(t,x) \circ W(\mathrm{d}x), \quad t \geq 0, x \in \mathbb{R};$$
$$u(0,x) = g(x), \quad x \in \mathbb{R}, \tag{6.1}$$

where the above stochastic integral should be understood in the sense of the anticipative Stratonovich integral (see [13] or [12]). However, since it is difficult to get any uniqueness result for such an equation, we prefer to rewrite it in a different form. In fact, since

$$\frac{1}{\sqrt{\varepsilon}}c\left(\frac{x}{\varepsilon}\right) = \overline{c}\frac{\partial W_\varepsilon}{\partial x}(x),$$

we can rewrite the original $u^\varepsilon$ equation as

$$\frac{\partial u^\varepsilon}{\partial t}(t,x) = \frac{1}{2}\frac{\partial}{\partial x}\left(a\left(\frac{\cdot}{\varepsilon}\right)\frac{\partial u^\varepsilon}{\partial x}\right)(t,x) + \overline{c}\frac{\partial}{\partial x}(W_\varepsilon u^\varepsilon)(t,x) - \overline{c}W_\varepsilon(x)\frac{\partial u^\varepsilon}{\partial x}(t,x),$$
$$u^\varepsilon(0,x) = g(x).$$

Hence we might expect that the limiting equation reads

$$\frac{\partial u}{\partial t}(t,x) = \frac{\tilde{a}}{2}\frac{\partial^2 u}{\partial x^2}(t,x) + \overline{c}\frac{\partial}{\partial x}(Wu)(t,x) - \overline{c}W(x)\frac{\partial u}{\partial x}(t,x), \quad t \geq 0, x \in \mathbb{R};$$
$$u(0,x) = g(x), \quad x \in \mathbb{R}. \tag{6.2}$$

Would the equation be posed on a compact interval, or would the random process $\{W(x), x \in \mathbb{R}\}$ have a.s. bounded trajectories, then the existence and the uniqueness for (6.2) would be very easy to obtain. Here we will content ourselves with the following result.

**Theorem 6.1.** *The parabolic PDE* (6.2) *has a solution* $u \in L^2_{\mathrm{loc}}(\mathbb{R}_+; H^1_{\mathrm{loc}}(\mathbb{R}))$ *a.s., which is given by the Feynman–Kac formula*

$$u(t,x) := \mathbb{E}\left[g(X_t^x)\exp\left(\frac{\overline{c}}{\tilde{a}}\int_\mathbb{R} L_t^{y-x}W(\mathrm{d}y)\right)\right]. \tag{6.3}$$

**Proof.** We need to show that $u$, given by the formula

$$u(t,x) = \mathbb{E}\left[g(X_t^x)\exp\left(\frac{\bar{c}}{\tilde{a}}\int_{\mathbb{R}} L_t^{y-x} W(\mathrm{d}y)\right)\right],$$

belongs to $L^2_{\mathrm{loc}}(\mathbb{R}_+; H^1_{\mathrm{loc}}(\mathbb{R}))$, and solves the parabolic PDE (6.2).

For that sake, we first define an approximation of the Wiener process $W$, by the formula

$$W_n(x) = [(W * \rho_n)(x) \wedge n] \vee (-n),$$

where $*$ stands for the convolution operation, and $\rho_n(x) := n\rho(nx)$, $\rho$ is a smooth map from $\mathbb{R}$ into $\mathbb{R}_+$ with compact support, such that

$$\int_{\mathbb{R}} \rho(x)\,\mathrm{d}x = 1.$$

Let $u^n$ denote the solution of the approximating PDE

$$\begin{aligned}
\frac{\partial u^n}{\partial t}(t,x) &= \frac{\tilde{a}}{2}\frac{\partial^2 u^n}{\partial x^2}(t,x) + \bar{c}W'_n(x)u^n(t,x), \quad t \geq 0, x \in \mathbb{R}; \\
u^n(0,x) &= g(x), \quad x \in \mathbb{R}.
\end{aligned} \quad (6.4)$$

It follows from the Feynman–Kac formula, see e.g. [9], that

$$\begin{aligned}
u^n(t,x) &= \mathbb{E}\left[g(x+X_t)\exp\left(\bar{c}\int_0^t \bar{W}'_n(x+X_s)\,\mathrm{d}s\right)\right] \\
&= \mathbb{E}\left[g(x+X_t)\exp\left(\frac{\bar{c}}{\tilde{a}}\int_{\mathbb{R}} \bar{W}'_n(z)L_t^{z-x}\,\mathrm{d}z\right)\right],
\end{aligned} \quad (6.5)$$

where $X_t = \sqrt{\tilde{a}}B_t$, and $\{B\}$ is a standard Brownian motion defined on $(\Omega, \mathcal{F}, \mathbb{P})$, while $L_s^y$ denotes its local time at time $s$ and point $y$. It follows from arguments similar to but simpler than those in Section 4 that

$$\mathbb{E}\left[g(x+X_t)\exp\left(\frac{\bar{c}}{\tilde{a}}\int_{\mathbb{R}} \bar{W}'_n(z)L_t^{z-x}\,\mathrm{d}z\right)\right] \to \mathbb{E}\left[g(x+X_t)\exp\left(\frac{\bar{c}}{\tilde{a}}\int_{\mathbb{R}} L_t^{z-x}W(\mathrm{d}z)\right)\right].$$

For each $M > 0$, we now write an equation satisfied by $u^n$:

$$\begin{aligned}
\frac{\partial u^n}{\partial t}(t,x) &= \frac{\tilde{a}}{2}\frac{\partial^2 u^n}{\partial x^2}(t,x) + \bar{c}\frac{\partial(W_n u^n)}{\partial x}(t,x) - \bar{c}W_n(x)\frac{\partial u^n}{\partial x}(t,x), \quad t > 0, -M < x < M \\
u^n(0,x) &= g(x), \qquad u^n(t,-M) = \xi^n(t,-M), \qquad u^n(t,M) = \xi^n(t,M),
\end{aligned} \quad (6.6)$$

where $\xi^n$ denotes the right-hand side of (6.5). It is now easy to show that

$$v^n(t,x) := u^n(t,x) - x\frac{\xi^n(t,M) - \xi^n(t,-M)}{2M} - \frac{\xi^n(t,M) + \xi^n(t,-M)}{2}$$

solves eq. (6.6) but with homogeneous Dirichlet boundary conditions. Now $v^n$ converges strongly in $L^2_{\mathrm{loc}}(\mathbb{R}_+; H^1_0(-M,M))$, $P$ a.s., towards the solution of the parabolic PDE

$$\begin{aligned}
\frac{\partial v}{\partial t}(t,x) &= \frac{\tilde{a}}{2}\frac{\partial^2 v}{\partial x^2}(t,x) + \bar{c}\frac{\partial(Wv)}{\partial x}(t,x) - \bar{c}W(x)\frac{\partial v}{\partial x}(t,x), \quad t \geq 0, -M < x < M; \\
v(0,x) &= g(x), \quad -M < x < M, \qquad v(t,-M) = v(t,M) = 0.
\end{aligned}$$

We conclude the $u := \lim_{n\to\infty} u^n$ belongs to the space $L^2_{\text{loc}}(\mathbb{R}_+; H^1_{\text{loc}}(\mathbb{R}))$ a.s., and it satisfies (6.2) in the variational sense, i.e. for any $t > 0$, any $\varphi \in C^2(\mathbb{R})$ with compact support, and a.s. ($\langle \cdot, \cdot \rangle$ denotes the scalar product in $L^2(\mathbb{R})$),

$$\langle u(t), \varphi \rangle = \langle g, \varphi \rangle + \int_0^t \left[ \frac{1}{2} \langle u(s), \varphi'' \rangle - \overline{c} \langle W u(s), \varphi' \rangle - \overline{c} \left\langle W \frac{\partial u}{\partial x}(s), \varphi \right\rangle \right] \mathrm{d}s. \qquad \Box$$

***Remark 6.2.*** *The same problem on a bounded interval $(a, b) \subset \mathbb{R}$ with Dirichlet homogeneous boundary conditions could be treated either by the same method, or by a much simpler PDE argument. Namely, noting that $\frac{1}{\sqrt{\varepsilon}} c(\frac{x}{\varepsilon}) = \overline{c} W'_\varepsilon(x)$, we have that $u^\varepsilon$ is the unique element of $L^2_{\text{loc}}(\mathbb{R}_+; H^1_0(a, b)) \cap C(\mathbb{R}_+; L^2(a, b))$ which satisfies the PDE*

$$\begin{cases} \frac{\partial u^\varepsilon}{\partial t} = \frac{1}{2} \frac{\partial}{\partial x} \left( a\left(\frac{x}{\varepsilon}\right) \frac{\partial u^\varepsilon}{\partial x} \right) + \overline{c} \frac{\partial}{\partial x}(W_\varepsilon u^\varepsilon) - \overline{c} W_\varepsilon \frac{\partial u^\varepsilon}{\partial x}, & a < x < b,\ t > 0; \\ u^\varepsilon(0, x) = g(x), \qquad u^\varepsilon(t, a) = u^\varepsilon(t, b) = 0. \end{cases}$$

*It is then not very difficult to show that $u^\varepsilon \Rightarrow u$, the unique solution again in $L^2_{\text{loc}}(\mathbb{R}_+; H^1_0(a, b)) \cap C(\mathbb{R}_+; L^2(a, b))$ of eq. (6.2) in $\mathbb{R}_+ \times (a, b)$, with homogeneous Dirichlet boundary conditions.*